\documentclass[12pt,twoside]{amsart}
\usepackage{amssymb}

\newtheorem{Thm}{Theorem}
\newtheorem{Cor}{Corollary}
\newtheorem{Lemma}{Lemma}
\newtheorem{Prop}{Proposition}

\theoremstyle{definition}
\newtheorem{Defn}{Definition}

\newcommand{\norm}[1]{\left\Vert#1\right\Vert}
\newcommand{\abs}[1]{\left\vert#1\right\vert}

\newcommand{\R}{\ensuremath{\mathbb R}}
\newcommand{\N}{\ensuremath{\mathbb N}}
\newcommand{\Z}{\ensuremath{\mathbb Z}}

\DeclareMathOperator{\Var}{Var}

\newcommand{\set}[1]{\left\{#1\right\}}

\newcommand{\Part}{\mathcal{P}}
\renewcommand{\phi}{\varphi}
\newcommand{\Prod}{\mathrm{Pr}}
\newcommand{\St}{\mathrm{St}}
\newcommand{\NC}{\mathit{NC}}
\newcommand{\Int}{\mathit{Int}}

\title[Free stochastic measures]{Free stochastic measures via noncrossing partitions}
\author[M. Anshelevich]{Michael Anshelevich}
\thanks{This work is supported in part by the Fannie and John Hertz
  Foundation Fellowship}
\address{Department of Mathematics, University at California \\
Berkeley, CA 94720, USA}
\email{mashel@math.berkeley.edu}
\subjclass{Primary 46L50; Secondary 60J30, 60G}

\date{October 27, 1999}

\begin{document}

\begin{abstract}
We consider free multiple stochastic measures in the combinatorial framework of the lattice
of all diagonals of an $n$-dimensional space. In this free case, one can restrict the
analysis to only the noncrossing diagonals. We give definitions of what free multiple
stochastic measures are, and calculate them for the free Poisson and free compound Poisson
processes. We also derive general combinatorial It\^{o}-type relationships between free
stochastic measures of different orders. These allow us to calculate, for example, free
Poisson-Charlier polynomials, which are the orthogonal polynomials with respect to the free
Poisson measure.
\end{abstract}

\maketitle

\section{Introduction}

The motivation for this paper is twofold. On the one hand, in \cite{RW97} Rota and
Wallstrom show that much of the classical theory of multiple stochastic integrals can be
done combinatorially, using the properties of the lattice of all partitions of a set,
especially the M\"{o}bius inversion formula. As one consequence they get a number of
combinatorial formulas describing the properties of orthogonal polynomials. On the other
hand, recently Biane and Speicher in \cite{BS98} made major advances in the study of the
free Brownian motion, started earlier in \cite{Bia97a,KS92,Spe91,Fag91}.

We continue the study of more general free stochastic processes, concentrating especially
on the free Poisson process. The starting point of the \cite{RW97} paper is the observation
(which has been made before) that the first difficulty in dealing with stochastic measures,
as compared with scalar measures, is that multiple stochastic measures should not be taken
as simply product measures of one-dimensional ones. Indeed, various diagonals in the
$n$-dimensional space, which have Lebesgue measure $0$, have nonzero product stochastic
measure. Rota and Wallstrom point out, however, that removal of these diagonals from the
space, which is how one usually defines multiple stochastic measures, corresponds precisely
to the M\"{o}bius inversion on the lattice of all partitions of a set, which is the same as
the lattice of all diagonals (see Section \ref{sec:diag}). We apply this idea to the study
of stochastic measures in free probability. Inspired by \cite{RW97}, we define
$n$-dimensional free stochastic measures and, more generally, free stochastic measures
depending on a partition $\pi$ in Section \ref{sec:Defn}. Interestingly, in the free case
the diagonals corresponding to crossing partitions all have weight $0$ to begin with. This
is in accordance with general approach of Speicher that combinatorially, the transition
from the classical to the free probability corresponds to the transition from the lattice
of all to the lattice of noncrossing partitions.

After the free Brownian motion, the most important process with free increments is the free
Poisson process. Using the combinatorial machinery, we can calculate explicitly the
multiple stochastic measures for the free Poisson process, and more generally for the free
compound Poisson processes. These in turn give us recurrence relations for the orthogonal
polynomials with respect to the corresponding (scalar) measures. In particular we calculate
the free Poisson-Charlier polynomials. Finally, in the language of the lattice of
noncrossing partitions one can easily express the It\^{o} product formula for general free
stochastic measures.

The paper is organized as follows. In Section \ref{sec:NotDef} we collect various
combinatorial preliminaries, and the main definitions. In Section \ref{sec:prel}, we look
at some general properties of free stochastic measures, and calculate the distributions for
their main diagonal measures. These diagonal measures are calculated explicitly for the
free Poisson process in Section \ref{sec:Poi} and for the free compound Poisson processes
in Section \ref{sec:cPoi}. Section \ref{sec:Prod} is devoted to the consideration of
product measures, especially in the free Brownian motion and free Poisson cases, and
combinatorial formulas that have implications for families of orthogonal polynomials. It
also contains the combinatorial It\^{o} product formula. In section \ref{sec:wd} we show
that at least for a particular scheme, the stochastic measures are always well-defined.
Section \ref{sec:KS} is devoted to various recursion relations between stochastic measures,
and the relation to orthogonal polynomials. Finally, in Section \ref{sec:func} we list a
few preliminary facts about the result of integration with respect to free stochastic
processes.

\noindent \textbf{Acknowledgments:}
I'd like to thank my advisor, Professor Dan-Virgil Voiculescu, for suggesting the free
Poisson process as a subject of investigation, for many helpful questions and suggestions
during the preparation of this paper, and general support throughout. I also thank
Professor Steve Evans for bringing the paper \cite{RW97} to my attention, and Daniel
Markiewicz for useful comments on a preliminary version of this paper.

\section{Notation and definitions}
\label{sec:NotDef}

\subsection{Partitions.}
\label{sec:part}
We will consider the following three lattices of partitions. By $\Part(n)$ we'll denote the
lattice of all partitions of the set $\set{1, 2, \ldots, n}$. By $\NC(n)$ we'll denote the
lattice of noncrossing partitions \cite{Kre72}. These are the partitions with the property
that
\begin{equation*}
i < j < k, \; i \stackrel{\pi}{\sim} k , \; j \stackrel{\pi}{\sim} l,  \; i
\stackrel{\pi}{\not \sim} j \; \Rightarrow
\; i < l < k.
\end{equation*}
Finally, the third lattice, used mostly for notational convenience, is the lattice
$\Int(n)$ of interval partitions \cite{vWa73}. These are the partitions whose classes are
intervals, and $\Int(n)$ is isomorphic as a lattice to the lattice of subsets of a set of
$(n-1)$ elements.

There is a partial order $\leq$ on $\Part(n)$ which restricts to the other two lattices. We
denote the smallest element in that order by $\hat{0} = \set{(1), (2), \ldots, (n)}$, and
the largest one by $\hat{1} = \set{(1, 2, \ldots, n)}$. We denote the meet and the join in
the lattices by $\wedge$ and $\vee$, respectively.

We need the following operations on partitions.
For $\pi \in \Part(n)$, we define $\pi^{op} \in \Part(n)$ to be $\pi$ taken in the opposite
order, i.e.
\begin{equation*}
i  \stackrel{\pi^{op}}{\sim} j \; \Leftrightarrow \; (n-i+1) \stackrel{\pi}{\sim} (n-j+1)
.
\end{equation*}
We define $\pi^k \in \Part(nk)$, the $k$-thickening of $\pi$, by
\begin{equation*}
i \stackrel{\pi^k}{\sim} j \; \Leftrightarrow \; [(i-1)/k] + 1 \stackrel{\pi}{\sim}
[(j-1)/k] + 1 ,
\end{equation*}
where $[ \cdot ]$ denotes the integer part of a real number. In words, we expand each point
of the set $\set{1, 2, \ldots, n}$ into $k$ points, and require that those points lie
consecutively and in the same class of $\pi^k$. For $\pi \in \Part(n), \sigma \in
\Part(k)$, we define $\pi + \sigma \in \Part(n+k)$ by
\begin{equation*}
i \stackrel{\pi + \sigma}{\sim} j \Leftrightarrow ((i, j \leq n, i \stackrel{\pi}{\sim} j)
\text{ or } (i, j > n, (i-n) \stackrel{\sigma}{\sim} (j-n))).
\end{equation*}
We'll denote $m \pi := \pi + \pi + \ldots + \pi \ $ $m$ times.

Following \cite{BLS96}, we divide the blocks of a noncrossing partition $\pi$ into inner
and outer: a block $B \in \pi$  is called \emph{inner} if there exist $(i
\stackrel{\pi}{\sim} j, \; i, j \not \in B)$ such that $i < k < j$ for some, hence all, $k
\in B$. A block that is not inner is called outer.

Finally, for $\pi \in \Part(n)$, we define the number of crossings of $\pi$, $c(\pi)$, to
be
\begin{equation*}
c(\pi) = \min (\abs{\sigma} - \abs{\pi} : \sigma \in \NC(n), \sigma \leq \pi).
\end{equation*}
In words, this is the minimal number of ``cuts'' in the classes of $\pi$ required to make
it noncrossing. Notice that this number is different from the reduced number of crossings
of \cite{Nic95}, the number of the restricted crossings of \cite{Bia97c}, and $m(\pi)$ of
\cite{Mar98}.

\subsection{Diagonals.}
\label{sec:diag}
Fix $n$, let $\pi$ be a set partition of $n$, with $|\pi| = k$ classes $B_1, B_2, \ldots,
B_k$. For a set $ S = \set{1, 2, \ldots, N}$, denote by $S^n_\pi$ the $\pi$-diagonal of
$S^n$, that is the set of $n$-tuples $(i_1, i_2, \ldots, i_n) \in S^n$ such that
\begin{equation*}
k \stackrel{\pi}{\sim} l  \ \text{(i.e. $k$ and $l$ lie in the same class of $\pi$)} \;
\Leftrightarrow \; i_k = i_l.
\end{equation*}
For example, to the partition of the set of eight elements $\set{(1, 5, 8), (2, 7), (3),
(4, 6)}$ there corresponds the set $\set{(i_1, i_2, i_3, i_4, i_1, i_4, i_2, i_1): i_1
\neq i_2 \neq i_3 \neq i_4}$.
Note that $S^n_\pi \neq \emptyset$ only if $N \geq \abs{\pi}$; more generally,
$\abs{\set{1, 2, \ldots, N}^k_\pi} = (N)_{\abs{\pi}}$, where $(N)_m = N(N-1) \cdots
(N-m+1)$. Similarly, denote by $S^n_{\geq \pi}$ the set of $n$-tuples $(i_1, i_2, \ldots,
i_n) \in S^n$ such that
\begin{equation*}
k \stackrel{\pi}{\sim} l  \Rightarrow i_k = i_l.
\end{equation*}
Note that
\begin{equation}
\label{gpi}
S^n_{\geq \pi} = \bigcup_{\sigma \geq \pi} S^n_\sigma
\end{equation}

\subsection{Processes with free increments.}
\label{sec:Defn}
Let $(\mathcal{A}, \phi)$ be a tracial $W^\ast$- noncommutative probability space. That is,
$\mathcal{A}$ is a finite von Neumann algebra, and $\phi$ is a faithful normal trace state
on it. We will call the elements of $\mathcal{A}$ noncommutative random variables, or
random variables for short.

For definitions of free probabilistic notions that are not given here we refer the reader,
for example, to the monograph \cite{VDN92}. Let $\mu$ be a freely infinitely divisible
distribution with compact support; normalize it so that $\Var (\mu) = 1$. Let
$\set{\mu_t}_{t \in [0, \infty)}$ be the corresponding additive free convolution semigroup.
Note that in the free case, as opposed to the classical case, both the free normal
distribution (the semicircular distribution) and the free Poisson distribution have compact
support, so the condition is not as restrictive as it might appear. We will return to the
matter of extending the contents of this paper to general freely infinitely divisible
distributions elsewhere.
\begin{Defn}
\label{defn:fi}
A stationary stochastic process with freely independent increments is a map from the set of
finite half-open intervals $I = [a, b) \subset \R$ to the self-adjoint part of
$(\mathcal{A},
\phi)$ (which can be extended to the map on all Borel subsets) $I \mapsto X_I$ with the
following three properties:
\begin{enumerate}
\item $I \cap J = \emptyset \Rightarrow$ $X_I$ and $X_J$ are freely independent
 (free increments),
\item $I_1 \cap I_2 = \emptyset, I_1 \cup I_2 = J \Rightarrow X_{I_1} + X_{I_2} = X_J$
 (additivity),
\item The distribution of $X_I$ is $\mu_{\abs{I}}$ (stationarity).
\end{enumerate}
\end{Defn}
Substantial study of processes with free increments was started in \cite{GSS92} and
\cite{Bia98}. In particular, it was shown in \cite{GSS92} that for any $\mu$ as above,
there is a realization of such a process. Note that throughout the paper the terms ``free
stochastic process'' and ``free stochastic measure'' will be used interchangeably.

We now define the product measures $\Prod_\pi(A)$ and the stochastic measures $\St_\pi(A)$,
depending on the partition $\pi$ of $k$. These will again be additive processes on the real
line with free identically distributed increments. Note that this is in contrast with
\cite{RW97}, where the corresponding objects are processes on a $k$-dimensional space.
\begin{Defn}
\label{defn:pr}
Let $A$ be a union of half-open intervals in $\R$. Denote $X = X_A$, and for an arbitrary
$N$ let $X_1^{(N)}, X_2^{(N)}, \ldots, X_N^{(N)}$ be freely independent, identically
distributed, and add up to $X$. Note that henceforth we will usually omit the explicit
dependence on $N$, to simplify notation. Then
\begin{equation*}
\St_\pi(A) = \lim_{N \rightarrow \infty} \sum_{\substack{(i_1, i_2, \ldots, i_k) \\
\in \set{1, 2, \ldots, N}_\pi^k}} X_{i_1}^{(N)} X_{i_2}^{(N)} \cdots X_{i_k}^{(N)},
\end{equation*}
\begin{equation*}
\Prod_\pi(A) = \lim_{N \rightarrow \infty} \sum_{\substack{(i_1, i_2, \ldots, i_k) \\
\in \set{1, 2, \ldots, N}_{\geq \pi}^k}} X_{i_1}^{(N)} X_{i_2}^{(N)} \cdots X_{i_k}^{(N)},
\end{equation*}
where the limit here and, unless noted otherwise, elsewhere are taken in the operator norm.
For the discussion of the existence of the limits, see Sections \ref{sec:prel} and
\ref{sec:wd}.
\end{Defn}

We call $\psi_k := \St_{\hat{0}}$ \emph{the} stochastic measure of degree $k$. We also
define the $k$-th diagonal measure of the process by
\begin{equation*}
\Delta_k (A) = \St_{\hat{1}} = \Prod_{\hat{1}} = \lim_{N \rightarrow \infty}
\sum_{i=1}^N X_i^k.
\end{equation*}
Note that the second diagonal measure of the process $\Delta_2 (A)$ is frequently called
the quadratic variation of the process and denoted by $\langle X, X \rangle$.

Throughout most of the paper we will fix the set $A$ and write $X = X_A$, $\St_\pi := \St_\pi
(A)$, etc.

\subsection{Multidimensional $R$-transform and noncrossing cumulants.}
This section could have been taken directly from, say, \cite{NS96a} and is included for
completeness.

Given a family $x_1, x_2, \ldots, x_k$ in a noncommutative probability space $(\mathcal{A},
\phi)$, their joint distribution is the collection of their joint moments
\begin{equation*}
M(x_{i_1}, x_{i_2}, \ldots, x_{i_n}) = \phi(x_{i_1} x_{i_2} \cdots x_{i_n})
\end{equation*}
for $1 \leq i_j \leq k$, $1 \leq j \leq n$. For a partition $\pi \in \NC(n)$, we define
\begin{multline*}
M_{\pi} (x_{i_1}, x_{i_2}, \ldots, x_{i_n}) \\
= \prod_{B \in \pi} M(x_{i(j_1)}, x_{i(j_2)},
\ldots, x_{i(j_{\abs{B}})} : j_1 < j_2 < \ldots < j_{\abs{B}},
\set{j_1, j_2, \ldots, j_{\abs{B}}} = B).
\end{multline*}
We define the joint free cumulants of $(x_1, x_2, \ldots, x_k)$, which together comprise
the multidimensional $R$-transform, recursively by
\begin{equation*}
M_\pi(x_{i_1}, x_{i_2}, \ldots, x_{i_n}) = \sum_{\substack{\sigma \in
\NC(n) \\ \sigma \leq \pi}} R_\sigma (x_{i_1}, x_{i_2}, \ldots, x_{i_n})
\end{equation*}
where again
\begin{multline*}
R_{\sigma} (x_{i_1}, x_{i_2}, \ldots, x_{i_n}) \\
= \prod_{B \in \sigma} R(x_{i(j_1)}, x_{i(j_2)},
\ldots, x_{i(j_{\abs{B}})} : j_1 < j_2 < \ldots < j_{\abs{B}},
\set{j_1, j_2, \ldots, j_{\abs{B}}} = B).
\end{multline*}
The main property of the $R$-transform is that
\begin{equation*}
(i_j \stackrel{\pi}{\sim} i_l , \; x_{i_j} \text{ and } x_{i_l} \text{ freely independent
})
\Rightarrow R_\pi(x_{i_1}, x_{i_2}, \ldots, x_{i_n}) = 0.
\end{equation*}
Denote by $K$ the Kreweras complement map on $\NC(n)$ \cite{Kre72,NS96a}. This is a certain
bijection on the lattice $\NC(n)$ and it follows from the main property of the
$R$-transform that for $\set{x_1, x_2, \ldots, x_n}$ freely independent from $\set{y_1,
y_2, \ldots y_n}$
\begin{equation*}
\phi(x_1 y_1 x_2 y_2 \cdots x_n y_n) = \sum_{\pi \in \NC(n)} R_{K(\pi)} (x_1, x_2, \ldots x_n)
M_\pi (y_1, y_2, \ldots, y_n).
\end{equation*}

If $x_{i_1} = x_{i_2} = \ldots = x_{i_n} = x$, we write $R_\pi(x, x, \ldots, x)
:= R_\pi(x)$. We denote the individual moments and free cumulants by $m_n (x) =
M_{\hat{1}_n} (x)$, $r_n (x) = R_{\hat{1}_n} (x)$. Note that $r_1 (x) = m_1 (x) = \phi(x)$.

If $x$ has distribution $\mu$ which is freely infinitely divisible, and $y$ has
distribution $\mu_t$, then $r_n(y) = t r_n(x)$. Therefore $R_\pi (y) = t^{\abs{\pi}} R_\pi
(x)$.

Finally, for an algebra element $Z$, we denote by $Z^\mathrm{o}$ the centered version of
$Z$, $Z^\mathrm{o}
= Z - \phi(Z)$. Note that $r_1(Z^\mathrm{o}) = 0, r_n(Z^\mathrm{o}) = r_n (Z)$ for $n > 1$.

\section{Preliminaries}
\label{sec:prel}

As defined, the product and stochastic measures depend on the particular triangular array
$\set{X_i^{(N)}}_{i=1}^N, N \in \N$. We will show that the limits exist for a particular
choice of this array in Section \ref{sec:wd}. For now, we make a number of  observations
which are consequences of the free independence of the increments of the process, and which
will hold for any such array.

\begin{Lemma}
\label{lem:non-c}
\begin{equation*}
\lim_{N \rightarrow \infty} \phi \biggl( \sum_{\substack{(i_1, i_2, \ldots, i_k) \\
\in \set{1, 2, \ldots, N}_\pi^k}} X_{i_1} X_{i_2} \cdots X_{i_k}  \biggr) =
\begin{cases}
R_\pi (X) & \text{if $\pi$ is noncrossing} \\
0 & \text{if $\pi$ is crossing}.
\end{cases}
\end{equation*}
\end{Lemma}

\begin{proof}
\begin{align}
\phi \biggl(  \sum_{\substack{(i_1, i_2, \ldots, i_k) \\
\in \set{1, 2, \ldots, N}_\pi^k}} X_{i_1} X_{i_2} \cdots X_{i_k} \biggr)
 &= \sum_{\substack{\sigma \in \NC(k) \\ \sigma \leq \pi}}
\sum_{\substack{(i_1, i_2, \ldots, i_k) \\
\in \set{1, 2, \ldots, N}_\pi^k}} R_\sigma
\left( X_{i_1},  X_{i_2},  \ldots,  X_{i_k} \right) \notag \\
 &= (N)_{\abs{\pi}} \sum_{\substack{\sigma \in \NC(k) \\ \sigma \leq \pi}} N^{-\abs{\sigma}}
R_\sigma (X, X, \ldots, X) \label{nc}
\end{align}
since $r_j(X_i) = \frac{1}{N} r_j(X)$. If $\pi$ is noncrossing, then the limit, as $N
\rightarrow
\infty$, of the expression
\eqref{nc} is $R_\pi (X)$. On the other hand, assume that $\pi$ is crossing. As in
\cite{Bia98}, the number of elements of $\NC(k)$, which is a Catalan number, is less than
$4^k$ and for each $\sigma$, $\abs{R_\sigma (X)} \leq 4^k \norm{X}^k$. Thus the absolute
value of the expression \eqref{nc} is less than $4^{2k} \norm{X}^k N^{- c(\pi)}$. In
particular, it converges to $0$ as $N \rightarrow
\infty$.
\end{proof}

\begin{Thm}
\label{thm:non-c}
\begin{equation*}
\St_{\pi} = \lim_{N \rightarrow \infty} \sum_{\substack{(i_1, i_2, \ldots, i_k) \\
\in \set{1, 2, \ldots, N}_\pi^k}}
 X_{i_1} X_{i_2} \cdots X_{i_k} = 0
\end{equation*}
unless $\pi$ is noncrossing.
\end{Thm}

\begin{proof}
More generally,
\begin{multline}
\phi \biggl(  \biggl( \Bigl( \sum_{\substack{(i_1, i_2, \ldots, i_k) \\
\in \set{1, 2, \ldots, N}_\pi^k}} X_{i_1} X_{i_2} \cdots X_{i_k} \Bigr)
\Bigl( \sum_{\substack{(i_1, i_2,
\ldots, i_k) \\ \in \set{1, 2, \ldots, N}_\pi^k}} X_{i_1} X_{i_2} \cdots X_{i_k} \Bigr)^\ast
\biggr)^n \biggr) \\
= \phi \biggl( \sum_{\substack{\sigma \in \Part(2nk) \\ \sigma \geq n(\pi + \pi^{op}) \\
\sigma \wedge 2n \hat{1}_k =
n (\pi + \pi^{op})\;}} \sum_{\substack{(i_1, i_2, \ldots, i_{2kn}) \\ \in
\set{1, 2, \ldots, N}_{\sigma}^{2kn}}} X_{i_1} X_{i_2} \cdots X_{i_{2kn}} \biggr) \label{x-xs}.
\end{multline}

For $\sigma$ as in equation \eqref{x-xs}, $c(\sigma) \geq 2n c(\pi)$. Applying Lemma
\ref{lem:non-c} and using the estimate in the proof of that Lemma, we see that the above
expression \eqref{x-xs} is less than $4^{4nk} \norm{X}^{2nk} N^{- 2n c( \pi)}
d^{\abs{\pi}}_n$, where $d^m_n = \abs{\set{\sigma \in \Part(2 n m) : \sigma \wedge 2n
\hat{1}_m = 2n \hat{1}_m}}$. It was shown in Theorem 5.3.4 of \cite{BS98} that $\lim_{n
\rightarrow \infty} (d^m_n)^{1/2n} = (m + 1)$ (note that our use of $m$ and $n$ is the
opposite of theirs). Therefore
\begin{multline*}
\norm{\sum_{\substack{(i_1, i_2, \ldots, i_k) \\ \in \set{1, 2, \ldots, N}_\pi^k}}
X_{i_1} X_{i_2} \cdots X_{i_k}} \\
 = \lim_{n \rightarrow \infty} \left[ \phi \biggl(  \biggl( \Bigl(
 \sum_{\substack{(i_1, i_2, \ldots, i_k) \\
\in \set{1, 2, \ldots, N}_\pi^k}} X_{i_1} X_{i_2} \cdots X_{i_k} \Bigr)
\Bigl( \sum_{\substack{(i_1, i_2,
\ldots, i_k) \\ \in \set{1, 2, \ldots, N}_\pi^k}} X_{i_1} X_{i_2} \cdots X_{i_k} \Bigr)^\ast
\biggr)^n \biggr) \right]^{1/2n} \\
 \leq 4^{2k} \norm{X}^{k} (\abs{\pi} + 1) N^{-c(\pi)}.
\end{multline*}
In particular,
\begin{equation*}
\lim_{N \rightarrow \infty} \norm{\sum_{\substack{(i_1, i_2, \ldots, i_k) \\
\in \set{1, 2, \ldots, N}_\pi^k}} X_{i_1} X_{i_2} \cdots X_{i_k}} = 0
\end{equation*}
unless $\pi$ is noncrossing.
\end{proof}

The following is the analog of Proposition 1 from \cite{RW97} for the lattice of
noncrossing partitions. Note that Proposition for the lattice of all partitions follows
directly from Definition \ref{defn:pr} and so remains true as well.
\begin{Cor}
\label{cor:defn}
The measures $\Prod_\pi$ and $\St_\pi$ are related as follows: for $\pi \in \NC(k)$
\begin{align*}
\Prod_\pi &= \sum_{\substack{\sigma \in \NC(k) \\ \sigma \geq \pi}} \St_\sigma, \\
\St_\pi &= \sum_{\substack{\sigma \in \NC(k) \\ \sigma \geq \pi}} \mu(\pi, \sigma)
\Prod_\sigma,
\end{align*}
where $\mu(\pi, \sigma)$ is the M\"obius function on the lattice of noncrossing partitions
\cite{Kre72}.
\end{Cor}
\begin{proof}
By Definition \ref{defn:pr} and equation \eqref{gpi}, $\Prod_\pi = \sum_{\sigma \in
\Part(k),
\sigma
\geq \pi} \St_\sigma$. By the Theorem, for $\sigma \not \in \NC(k), \St_\sigma = 0$. The
second equality follows from the first one by the use of M\"{o}bius inversion on the
lattice of noncrossing partitions.
\end{proof}

Hereafter, it is reasonable to consider only product and stochastic measures
corresponding to noncrossing partitions, for the following reason. Following \cite{RW97},
we call a measure \emph{multiplicative} if for any partition $\pi$ (of the type to be
determined below)
\begin{equation}
\label{eq:mult}
\phi(\St_\pi) = \prod_{B \in \pi} \phi(\Delta_{\abs{B}}).
\end{equation}
It is easy to see that if the process is not commutative, the measures are not in general
multiplicative with respect to the lattice of all partitions. For example, let $\pi
= \set{(1,3), (2,4)}$. The partition is crossing, so the left-hand-side of equation
\eqref{eq:mult} is $0$, while for example for the free Brownian motion the right-hand-side
is $\phi(\Delta_2(A))^2 = \abs{A}^2 \neq 0$.


\begin{Cor}
\label{cor:mult}
Stochastic
processes with freely independent increments are multiplicative with respect to the lattice
of noncrossing partitions. That is, for $\pi \in \NC(k)$, $\phi(\St_\pi) =
\prod_{B \in \pi} \phi(\Delta_{\abs{B}})$.
\end{Cor}
\begin{proof}
By Lemma \ref{lem:non-c},
\begin{equation*}
\phi(\St_\pi) = R_\pi(X) = \prod_{B \in \pi} r_{\abs{B}}(X) = \prod_{B \in \pi} \phi(\Delta_{\abs{B}}).
\end{equation*}
Here the last equality follows from applying the first equality to $\pi = \hat{1}$.
\end{proof}

An immediate consequence
of Corollary \ref{cor:defn} is the analog of Theorem 1 from \cite{RW97}, which expresses
the stochastic measure of degree $k$ as a linear combination of various product measures.
In fact more is true, as in this particular case not even all noncrossing partitions
contribute to the sum. Call a \emph{singleton} a class of a partition consisting of one element.

\begin{Prop}
\label{prop:center}
For a
noncrossing partition $\pi$ that contains an inner singleton, $\St_\pi = 0$ if the process
is centered, that is, $\phi(X)=0$.
\end{Prop}

\begin{proof}
If $\pi$ contains a singleton, and the process is centered, then
\begin{equation*}
\phi \biggl( \sum_{\substack{(i_1, i_2, \ldots, i_k) \\
\in \set{1, 2, \ldots, N}_\pi^k}} X_{i_1} X_{i_2} \cdots X_{i_k}  \biggr) = 0.
\end{equation*}
In fact each term is $0$, for any $N$. On the other hand, if $\pi$ contains an inner
singleton, so does any partition $\sigma \geq n(\pi + \pi^{op}), \sigma \wedge (2n
\hat{1}_k)
= n(\pi + \pi^{op})$ with at most $(n-1)$ crossings, for any $n$. Thus in the sum in the
equation \eqref{x-xs}, only the partitions $\sigma$ with $c(\sigma) \geq n$ enter, and so
by the argument following that equation we see that the limit defining $\St_\pi$ is $0$.
\end{proof}

The following statement fits well with the original definition of free independence of Voiculescu
\cite{Voi85,VDN92}.
\begin{Cor}
If the process is centered,
\begin{equation*}
\psi_k = \St_{\hat{0}} = \lim_{N \rightarrow \infty} \sum_{\substack{i_1 \neq i_2 \neq
\cdots \neq i_k \\
\text{all distinct}}}^N X_{i_1} X_{i_2} \cdots X_{i_k} = \lim_{N \rightarrow \infty}
\sum_{\substack{i_1 \neq i_2 \neq \cdots \neq i_k \\ \text{neighbors distinct}}}^N X_{i_1}
X_{i_2} \cdots X_{i_k}.
\end{equation*}
\end{Cor}

\begin{proof}
The only noncrossing partition with no inner singletons for which no consecutive elements lie in the same class is
$\hat{0}$.
\end{proof}

Finally, note that for any triangular array as above, we always have convergence in
distribution  for the diagonal measures (and hence for the stochastic measures by the
results in Section \ref{sec:KS}):

\begin{Thm}
\label{thm:r_n}
The free cumulants of the $k$-th diagonal measure of the
process are given by
\begin{equation*}
r_n(\Delta_k) = \lim_{N \rightarrow \infty} r_n \bigl(\sum_{i=1}^N X_i^k \bigr) =
r_{nk}(X).
\end{equation*}
\end{Thm}

\begin{proof}
For $\mu(\pi) := \mu(\pi, \hat{1})$ the M\"obius function on the lattice of noncrossing
partitions,
\begin{align*}
r_n \bigl( \sum_{i=1}^N X_i^k \bigr)
 &= \sum_{i=1}^N \sum_{\pi \in \NC(n)} \mu(\pi) M_\pi (X_i^k)\\
 &= \sum_{i=1}^N \sum_{\pi \in \NC(n)} \mu(\pi) \prod_{B_j \in \pi} m_{\abs{B_j}} (X_i^k)
\\
 &= \sum_{i=1}^N \sum_{\pi \in \NC(n)} \mu(\pi) \prod_{B_j \in \pi} m_{k \abs{B_j}} (X_i) \\
 &= \sum_{i=1}^N \sum_{\pi \in \NC(n)} \mu(\pi) \prod_{B_j \in \pi} \sum_{\sigma_j \in
\NC(k \abs{B_j})} \prod_{A \in \sigma_j} r_{\abs{A}} (X_i) \\
 &= \sum_{\pi \in \NC(n)} \mu(\pi) \prod_{B_j \in \pi} \sum_{\sigma_j \in
\NC(k \abs{B_j})} N^{1 - \sum_{j=1}^{\abs{\pi}} \abs{\sigma_j}} \prod_{A \in \sigma_j}
r_{\abs{A}} (X) \\
 &= \mu(\hat{1}) r_{nk}(X) + O(1/N) \\
 &= r_{nk}(X) + O(1/N).
\end{align*}
\end{proof}

\section{Diagonal measures: examples}

\subsection{Free Brownian motion}
It follows from the results of \cite{BS98} that for the free Brownian motion the quadratic
variation is the scalar process $I \mapsto \abs{I}$, and the higher diagonal measures are
$0$.

\subsection{Free Poisson process}
\label{sec:Poi}

The free Poisson distribution is a distribution obtained by the Poisson-type limit process
for freely independent variables \cite{VDN92}. It is characterized by the
property that all its free cumulants are equal to $1$. For the free Poisson process, a
remarkable representation was given in \cite{NS96a} (see also the Appendix to that paper).
Let $I \mapsto p_I$ be a projection-valued process $[0,1] \rightarrow \mathcal{A}$. That
is, $\left( I \cap J = \emptyset \Rightarrow p_i \perp p_J \right)$, $\left( I_1 \cap I_2 =
\emptyset, I_1 \cup I_2 = J \Rightarrow p_{I_1} + p_{I_2} = p_J\right)$, $\phi(p_I) =
\abs{I}$. Let $s \in \mathcal{A}$ be an element with a standard semicircular distribution,
freely independent from the family $\set{p_I}$. Then $I \mapsto \set{s p_I s}$ is a free
Poisson process for $I \subset [0, 1]$. To get the full process, pick a countable
collection of such families of projections, $p_I^{(n)}, I \subset [0, 1], n \in \Z$, which
are freely independent from each other. Then $ [t_1, t_2) \mapsto s \left( p_{[t_1 - [t_1],
1)}^{([t_1])} + \sum_{i = [t_1] + 1}^{[t_2] - 1} p_{[0, 1)}^{(i)} + p_{[0, t_2 -
[t_2])}^{([t_2])} \right) s$ is a free Poisson process. Here $[\cdot]$ again denotes the
integer part of a real number.

By Theorem \ref{thm:r_n} we know the distributions of its diagonal measures. In fact in
this case we can prove convergence in norm to a specific limit. First we have a technical
theorem.

\begin{Thm}
\label{thm:eze}
Let $Z_1, Z_2, \ldots, Z_k$ be centered, and fix $e$, a self-adjoint element freely
independent from the family $\set{Z_i}$. For $N \in \N$, let $e_1^{(N)}, e_2^{(N)},
\ldots, e_N^{(N)}$ (where henceforth we will again omit the dependence on $N$) be
self-adjoint, freely independent from the family $\set{Z_i}$ and be an orthogonal family,
that is, $e_i e_j = 0$ for $i \neq j$. In addition, let the $e_i$'s be identically
distributed and $\sum_{i=1}^N e_i = e$. Then for arbitrary indices $m_1, m_2, \ldots,
m_{k+1} \in \N$
\begin{equation*}
\lim_{N \rightarrow \infty} \sum_{i=1}^N
 e_i^{m_1} Z_1 e_i^{m_2} Z_2 \cdots e_i^{m_k} Z_k e_i^{m_{k+1}} = 0.
\end{equation*}
\end{Thm}

\begin{proof}
First note that since the $e_i$'s are orthogonal and identically distributed,
\begin{equation*}
\phi(e^n) = \phi \left( ( \sum_{i=1}^N e_i ) ^n \right) =
 \phi \left( \sum_{i=1}^N e_i^n \right) = N \phi(e_1^n).
\end{equation*}
\begin{align}
& \phi \left( \left( \sum_{i=1}^N e_i^{m_1} Z_1 e_i^{m_2} Z_2 \cdots e_i^{m_k}
Z_k e_i^{m_{k+1}} \right)^n \right)
 = \phi \left( \sum_{i=1}^N \left( e_i^{m_1} Z_1 e_i^{m_2} Z_2 \cdots e_i^{m_k}
Z_k e_i^{m_{k+1}} \right)^n \right) \notag \\
 & \qquad = \phi \left( \sum_{i=1}^N \left( e_i^{m_1+ m_{k+1}} Z_1 e_i^{m_2} Z_2 \cdots
e_i^{m_k} Z_k \right)^n \right) \notag \\
 & \qquad = \sum_{i=1}^N \sum_{\pi \in \NC(nk)} R_{K(\pi)} (Z_1, Z_2, \ldots, Z_k, Z_1, \ldots, Z_k)
\label{eze} \\
 & \qquad \qquad \qquad \qquad \times M_{\pi} (e_i^{m_1+ m_{k+1}}, e_i^{m_2}, \ldots,
 e_i^{m_k}, e_i^{m_1+ m_{k+1}}, \ldots, e_i^{m_k}) \notag.
\end{align}
$Z_i$'s are centered, so $r_1 (Z_i) = 0$. Thus only those partitions $\pi$ contribute to
the sum \eqref{eze} for which $K(\pi)$ has no single-element classes. In particular
$\abs{K(\pi)} \leq (nk)/2$ and so $\abs{\pi} \geq (nk)/2$ (since $\abs{K(\pi)} +
\abs{\pi} = nk+1$). Then the sum
\eqref{eze} is
\begin{multline*}
\sum_{i=1}^N \sum_{\substack{\pi \in \NC(nk) \\ \abs{\pi} \geq (nk)/2}}
R_{K(\pi)} (Z_1, Z_2, \ldots, Z_k, Z_1, \ldots, Z_k) N^{- \abs{\pi}} \\
 \times M_{(m_1 + m_{k+1}, m_2,
\ldots, m_k, m_1 + m_{k+1}, \ldots, m_k) \pi} (e, e, \ldots, e),
\end{multline*}
where $(u_1, u_2, \ldots, u_{nk}) \pi \in \NC(\sum_{i=1}^{nk} u_i)$ is the partition
obtained by expanding the $i$-th point of the set on which $\pi$ operates into $u_i$
points. That is,
\begin{multline*}
(i_1 \stackrel{(u_1, u_2, \ldots, u_{nk}) \pi}{\sim} i_2) \Leftrightarrow \\
(\sum_{j=1}^{s_1} u_j + 1
\leq i_1 \leq \sum_{j=1}^{s_1 + 1} u_j, \qquad\sum_{j=1}^{s_2} u_j + 1 \leq i_2 \leq
\sum_{j=1}^{s_2 + 1} u_j, \text{ and } s_1 \stackrel{\pi}{\sim} s_2).
\end{multline*}
(Note that $(m, m, \ldots, m) \pi = \pi^m$.) Therefore the sum \eqref{eze} is less in
absolute value than
\begin{equation*}
4^{2nk} (\max \norm{Z_i})^{nk} \norm{e}^{n \sum_{i=1}^{k+1} m_i} N^{1 - nk/2}.
\end{equation*}
We can always choose $Z_i$'s to be self-adjoint, while the $e_i$'s are so by assumption.
Hence
\begin{equation*}
\norm{\sum_{i=1}^N e_i^{m_1} Z_1 e_i^{m_2} Z_2 \cdots e_i^{m_k} Z_k e_i^{m_{k+1}}}
\leq 4^{2k} (\max \norm{Z_i})^k \norm{e}^{\sum_{i=1}^{k+1} m_i} N^{-k/2}
\end{equation*}
and so converges to $0$ as $N \rightarrow \infty$.
\end{proof}

\begin{Cor}
\label{cor:poi-d_n}
For the free Poisson process, the $k$-th diagonal measure is equal to the process itself.

\end{Cor}
\begin{proof}
First consider the process on the interval $[0, 1]$. In the notation from the beginning of
the section, $X = sps$ and $X_i = s p_i s$.
\begin{align*}
\sum_{i=1}^N (p_i s^2 p_i)^k
 &= \sum_{i=1}^N p_i s^2 p_i s^2 \cdots p_i s^2 p_i \\
 &= \sum_{i=1}^N \sum_{\pi \in \Int(k+1)} \phi(s^2)^{k- \abs{\pi} +1} p_i ((s^2)^\mathrm{o}
p_i)^{\abs{\pi} - 1} \\
 &= \sum_{i=1}^N \sum_{j=0}^k \binom{k}{j} p_i ((s^2)^\mathrm{o} p_i)^j
\end{align*}
since $\phi(s^2) = 1$. In the limit $N \rightarrow \infty$, by the Theorem the only term
that survives is the one for $j=0$ (i.e. $\pi = \hat{1}, \abs{\hat{1}} = 1$), and that term
is $\sum_{i=1}^N p_i
= p$.

Therefore
\begin{equation*}
\lim_{N \rightarrow \infty} \sum_{i=1}^N (sp_i s)^k = s \lim_{N \rightarrow \infty}
\sum_{i=1}^N (p_i s^2 p_i)^{k-1} s = sps.
\end{equation*}
For the full Poisson process, the result follows from the free independence of the
increments corresponding to disjoint intervals.
\end{proof}

In fact we can calculate explicitly all the product measures for the free Poisson process;
see Section \ref{sec:Prod}.

\subsection{Free compound Poisson processes}
\label{sec:cPoi}

Let $e_I$ be a process with identically distributed orthogonal increments on $[0, 1]$. It
is a map from all finite half-open intervals $I \subset [0, 1]$ to the self-adjoint part of
$(\mathcal{A}, \phi)$ such that $\left( I \cap J = \emptyset \Rightarrow e_I e_J = 0
\right)$, $(I_1 \cap I_2 = \emptyset, I_1 \cup I_2 = J \Rightarrow e_{I_1} + e_{I_2}
= e_J)$, and the distribution of $X_I$ depends only on $\abs{I}$. Note that this
implies that if $\mu$ is the distribution of  $e_{[0, 1)}$ then the distribution of $e_I$
is $\left( (1 - \abs{I}) \delta_0 + \abs{I} \mu \right)$. Let $s$ be a random variable with
the  standard semicircle distribution, freely independent from this process. Then by
\cite{NS96a} the process $s e_I s$ is a process with free increments. Its distribution is
\begin{equation*}
r_n(s^2 e)
 = \sum_{\pi \in \NC(n)} R_\pi (e) R_{K(\pi)} (s^2)
 = \sum_{\pi \in \NC(n)} R_\pi (e)
 = m_n(e).
\end{equation*}

That is,  the process $s e_I s$ is a free compound Poisson process (these have been
considered previously in \cite{GSS92} and \cite{Spe98}). In fact, that will be true even if
the process does not have moments that are finite. Indeed, for general distributions we can
use the ``free Fourier transform'', which has appeared for example in \cite{NS97a}. Namely,
there is a relationship between Voiculescu's $R$- and $S$-transforms: denoting $\alpha(z) =
z R(z)$, one has $S(w) = w^{-1} \alpha^{-1}(w)$.

\begin{Lemma}
\label{lem:S-tr}
Let $x$ be a random variable, and $s$ a standard semicircular random variable freely
independent from $x$. Let $y = sxs$. Then $S_y (w) = \frac{1}{w+1} S_x(w)$. Therefore
$R_y(w) = w^{-1} \psi_x (w)$. Consequently the distribution of $y$ is a free compound
Poisson distribution with the Levy measure equal to the distribution of $x$.
\end{Lemma}

\begin{proof}
For a free Poisson element $s^2$,  $R(z) = \frac{1}{1-z}$, so $S(w) = \frac{1}{1+w}$.
Therefore $S_y (w) = \frac{1}{1+w} S_x(w)$ and so $\alpha_y (w) = \frac{w}{1+w} S_x(w) =
\chi_x(w)$, which implies $R_y (z) = z^{-1} \psi_x (z)$.
\end{proof}

For a free compound Poisson process given in the standard form $s e_I s$, we
will call the process with orthogonal increments $e_I$ the \emph{generator} of the process
(cf. \cite{GSS92}).

To extend the process from $[0, 1]$ to the whole real line we can, once again, take a
countable family of processes $e_I^{(n)}$ with orthogonal increments on $[0, 1]$ which are
freely independent from each other and, thinking of $e_I^{(n)}$ as acting on $[n, n+1]$,
take their sum.

Again we know the distributions of the diagonal measures from Theorem \ref{thm:r_n}, but
can in fact find them explicitly:

\begin{Cor}
\label{cor:cPoi}
The diagonal measures of a free compound Poisson process with generator $e$ are free compound
Poisson processes with generators $e^k$.
\end{Cor}

\begin{proof}
\begin{align*}
\sum_{i=1}^N (s e_i s)^k
 &= s \sum_{i=1}^N e_i s^2 e_i s^2 \cdots e_i s \\
 &= s \sum_{i=1}^N \sum_{\pi \in \Int(k)} \phi(s^2)^{k - \abs{\pi}}
e_i^{\abs{B_1}} (s^2)^\mathrm{o} e_i^{\abs{B_2}} (s^2)^\mathrm{o} \cdots (s^2)^\mathrm{o}
e_i^{\abs{B_{\abs{\pi}}}} s,
\end{align*}
where $B_1, B_2, \ldots, B_{\abs{\pi}}$ are the blocks of $\pi$. By Theorem \ref{thm:eze},
once again the only term that survives in the limit $N \rightarrow \infty$ is the one for
$\pi = \hat{1}, \abs{\hat{1}} = 1$, and that term is $s \sum_{i=1}^N e_i^k s = s
(\sum_{i=1}^N e_i)^k s = s e^k s$.
\end{proof}

\section{Product measures}
\label{sec:Prod}

The key point of the paper \cite{RW97} is that the stochastic measures $\St_\pi$ can be
expressed combinatorially through the product measures $\Prod_\sigma$, and the latter are
indeed product vector measures. In the free case the situation is more complicated.

For a noncrossing partition $\pi$, let $i(\pi)$ and $o(\pi)$ be, respectively, the numbers
of the inner and outer classes of $\pi$; in particular, $i(\pi) + o(\pi) =
\abs{\pi}$. Here's another description of $i(\pi)$ and $o(\pi)$. There is a partial order on the
classes of the partition $\pi$, by the height, i.e. for two classes $B
\neq C$ of $\pi$, $B > C \Leftrightarrow \exists i, j \in B, k \in C : i < k < j$. Then
$i(\pi)$ is the number of edges in the incidence graph of this partial order, and $o(\pi)$
is the number of maximal elements under this order.

\begin{Prop}
For $X$ the free Brownian motion, for a partition $\pi$ with no inner singletons
\begin{equation*}
\Prod_\pi (A) = X^{\abs{\set{B \in \pi: \abs{B} = 1}}} \abs{A}^{\abs{ \set{B \in \pi: \abs{B}
 = 2}}} 0^{\abs{ \set{B \in \pi: \abs{B} > 2}}} = \prod_{B \in \pi} \Delta_{\abs{B}} (A).
\end{equation*}
That is, $\Prod_\pi$ are indeed product measures.
\end{Prop}

\begin{proof}
The proof will be by induction on the level in the above partial order on the classes of
the partition $\pi$. For the product measure $\Prod_\pi = \lim_{N \rightarrow
\infty} \sum_{(i_1, i_2, \ldots, i_k) \in \set{1, 2, \ldots, N}_{\geq \pi}^k} X_{i_1}
X_{i_2} \cdots X_{i_k}$ one can sum over the indices corresponding to the minimal classes,
which are precisely the classes that are intervals, without disturbing the rest of the
partition. If any of these classes contains at least $3$ elements, the result is $0$, while
each two-element class gives $\abs{A}$, which is a scalar and can be factored out. In the
partition resulting from factoring out the minimal two-element classes, all the classes
which were at the distance of $1$ from the minimum now become minimal.
\end{proof}

For a partition $\pi$ \emph{with} an inner singleton, $\Pr_\pi (A)= \sum_{\sigma \geq
\pi} \St_\sigma (A)$. Each of such $\sigma$'s will either contain an inner singleton, or a
class of at least 3 elements. In either case $\St_\sigma = 0$. Thus in this case $\Prod_\pi
= 0$, and so need not in general be a product measure.

\begin{Prop}
\label{prop:Poi}
Let $X$ be the free Poisson process. Suppose $\pi$ is a partition with the following property: if $U$, $V$ are inner classes of $\pi$ covered by a class $W$ and such that $\forall u \in U, v \in V: u < v$, then there exists $w \in W$ such that $\forall u \in U, v \in V: u < w < v$. Then 
\begin{equation*}
\Prod_\pi (A) = X^{o(\pi)} (1 + \abs{A})^{i(\pi)}.
\end{equation*}
In particular, $\Prod_\pi$ in this case are not product measures. Note, however, that if
$\pi$ is an interval partition, then $i(\pi) = 0$ and $o(\pi) = \abs{\pi}$, and so one does
again get a product measure decomposition.
\end{Prop}

\begin{proof}
Once again, the proof is by induction on the level of the class of $\pi$ in the partial
order. Since the diagonal measures of the free Poisson process are equal to the process
itself, each minimal (interval) class can be shrunk to a one-element class. Again, consider
a class which used to be at a distance of $1$ from the minimum, and after the above
shrinking, covers only a number of singletons:
\begin{multline*}
\sum_{i, j_1, j_2, \ldots, j_k, \ldots = 1}^N \cdots X_i X_{j_1} X_i X_{j_2} X_i
 \cdots X_{j_k} X_i \cdots \\
= \sum_{\substack{\sigma \in \Int(k+1) \\ \sigma = (B_1, B_2, \ldots, B_{\abs{\sigma}} ) }}
 \sum_{\substack{ i, j_1, j_2, \ldots, j_{\abs{\sigma} - 1} , \ldots = 1 \\
i \neq j_1, \ldots, i \neq j_{\abs{\sigma} - 1}}}^N \cdots X_i^{\abs{B_1}} X_{j_1}
X_i^{\abs{B_2}} \cdots  X_{j_{\abs{\sigma} - 1}} X_i^{\abs{B_{\abs{\sigma}}}} \cdots.
\end{multline*}
Even though the free Poisson process is not centered, it is easy to see that each inner
singleton contributes a factor of $\phi(X) = \abs{A}$. Thus all the singletons covered by
the $i$ class contribute $\sum_{\sigma \in \Int(k+1)} \abs{A}^{\abs{\sigma} - 1} = (\abs{A}
+ 1)^k$. Note that $k$ is precisely the number of classes covered by the $i$ class, i.e.
the number of edges emanating down from it in the partial order incidence graph. The result
easily follows by induction.
\end{proof}

Now we consider products of free stochastic measures.

\begin{Defn}
For $\pi \in \Part(k)$, define
\begin{equation*}
\overrightarrow{\prod}_{B \in \pi} \psi_{\abs{B}} =
\lim_{N \rightarrow \infty} \sum X_{i_1} X_{i_2} \cdots X_{i_k},
\end{equation*}
where the sum is taken over all collections of indices in $\set{1, 2, \ldots, N}^k$ such that
\begin{equation}
\label{ind}
m \stackrel{\pi}{\sim} l  \Rightarrow i_m \neq i_l.
\end{equation}
That is, on each of the classes of $\pi$ we consider the measure $\St_{\hat{0}}$. In
particular, if $\pi$ is an interval partition, then the result is the ordered product of
measures corresponding to all the classes.
\end{Defn}
The following is the analog of Theorem 4 of \cite{RW97}. This is a combinatorial form of
the general It\^{o} product formula.

\begin{Prop}
\label{prop:prod}
Let $\pi \in \Part(k)$ Then
\begin{equation}
\label{wedge}
\overrightarrow{\prod}_{B \in \pi} \psi_{\abs{B}} =
\sum_{\substack{\sigma \in \Part(k) \\ \sigma \wedge \pi = \hat{0}}} \St_\sigma =
\sum_{\substack{\sigma \in \NC(k) \\ \sigma \wedge \pi = \hat{0}}} \St_\sigma.
\end{equation}
\end{Prop}

\begin{proof}
Any $k$-tuple of indices determines a partition $\sigma \in \Part(k)$ by $(m
\stackrel{\sigma}{\sim} l)
\Leftrightarrow (i_m = i_l)$. Condition \eqref{ind} restricts $\sigma$ precisely to those
for which $\sigma \wedge \pi = \hat{0}$, and the union of all collections of indices
corresponding to such $\sigma$'s is precisely the set given by condition \eqref{ind}. This
decomposition of the set of indices gives the representation \eqref{wedge}.

The second equality follows from Theorem \ref{thm:non-c}.
\end{proof}

\begin{Cor}
For $\pi \in \Part(k)$,
\begin{equation*}
\overrightarrow{\prod}_{B \in \pi} \psi_{\abs{B}} =
\sum_{\hat{0} \leq \sigma \leq \pi} \mu_{\Part}(\hat{0}, \sigma) \Prod_\sigma =
\sum_{\substack{\sigma \in \NC(k) \\ \hat{0} \leq \sigma \leq \pi}}
\mu_{\NC}(\hat{0}, \sigma) \Prod_\sigma.
\end{equation*}
where $\mu_{\Part}$ and $\mu_{\NC}$ are the M\"{o}bius functions on the lattices of all and
of noncrossing partitions, respectively.
\end{Cor}
\begin{proof}
The proof of Theorem 4 in \cite{RW97} shows that
\begin{equation*}
\sum_{\hat{0} \leq \sigma \leq \pi} \mu(\hat{0}, \sigma) \Prod_\sigma =
\sum_{\substack{\sigma \in \Part(k) \\ \sigma \wedge \pi = \hat{0}}} \St_\sigma.
\end{equation*}
That proof is purely combinatorial, and works for noncommutative stochastic measures as
well as for the commutative ones. Therefore the first equality follows from Proposition
\ref{prop:prod} above. Moreover, the proof also works for the lattice of noncrossing
partitions, provided we use the appropriate M\"{o}bius function. Therefore the second
equality holds as well.
\end{proof}

\section{Orthogonality relations and the existence of the limits}
\subsection{Orthogonality relations}
The following Proposition is the analog of Proposition 9, Theorem 9, and Proposition 10
of \cite{RW97}.

\begin{Prop}
\label{prop:orth}
\

\begin{enumerate}
\item Let $\pi \in \NC(n)$. Then
\begin{equation*}
\phi \biggl( \overrightarrow{\prod}_{B \in \pi} \psi_{\abs{B}} \biggr) =
\sum_{\substack{\sigma \in \NC(n) \\
 \sigma \wedge \pi = \hat{0}}} \prod_{B \in \sigma} \phi(\Delta_{\abs{B}}).
\end{equation*}
\item Assume that the process is centered. Then
\begin{equation*}
\phi(\psi_n \psi_m) =
\begin{cases}
\phi(\Delta_2)^n = r_2 (X)^n = \abs{A}^n, \text{ if } m=n \\
0, \text{ if } m \neq n.
\end{cases}
\end{equation*}
\end{enumerate}
\end{Prop}

\begin{proof}
The first part is obtained by taking expectations of both sides in Proposition
\ref{prop:prod} and applying the multiplicativity property of Corollary \ref{cor:mult}.
The second part follows from the first part.
\end{proof}

\subsection{The existence of the limits}
\label{sec:wd}

By the results in the previous subsection, standard It\^{o} isometry arguments show that
the limit defining $\psi_n$ exists in $L^2$. But in fact, we want to show that the limits
exist, in norm, for all stochastic measures.

Without loss of generality, let $A = [0, 1)$. Choose two different partitions of the
interval $A$: $X_i = X_{[\frac{i-1}{N}, \frac{i}{N})}, i = 1, 2, \ldots, N$, $Y_j =
X_{[\frac{j-1}{M}, \frac{j}{M})}, j = 1, 2, \ldots, M$. Let  $Z_m = X_{[\frac{m-1}{MN},
\frac{m}{MN})}, m = 1, 2, \ldots, MN$ be their common refinement. Fix a partition $\pi$.
\begin{multline}
\left[\phi \biggl( \biggl( \sum_{\substack{(i_1, i_2, \ldots, i_k) \\
\in \set{1, 2, \ldots, N}_\pi^k}} X_{i_1} X_{i_2} \cdots X_{i_k} - \sum_{\substack
{(j_1, j_2, \ldots, j_k) \\
\in \set{1, 2, \ldots, M}_\pi^k}} Y_{j_1} Y_{j_2} \cdots Y_{j_k} \biggr)^n \biggr)
\right]^{1/n} \\
= \left[\phi \biggl( \biggl( \sum_{\substack{(i_1, i_2, \ldots, i_k) \\
\in \set{1, 2, \ldots, N}_\pi^k}}
(\sum_{s_1 =(i_1 - 1)M + 1}^{i_1 M} Z_{s_1}) (\sum_{s_2 =(i_2 - 1)M + 1}^{i_2 M} Z_{s_2})
\cdots (\sum_{s_k =(i_k - 1)M + 1}^{i_k M} Z_{s_k}) \right. \label{Cauchy} \\
- \left. \sum_{\substack{(j_1, j_2, \ldots, j_k) \\ \in \set{1, 2, \ldots, M}_\pi^k}}
(\sum_{t_1 =(j_1 - 1)N + 1}^{j_1 N} Z_{t_1}) (\sum_{t_2 =(j_2 - 1)N + 1}^{j_2 N} Z_{t_2})
\cdots (\sum_{t_k =(j_k - 1)N + 1}^{j_k N} Z_{t_k}) \biggr)^n \biggr) \right]^{1/n}.
\end{multline}
This expression can be written as a sum over partitions $\sigma \in \Part(nk)$ of free
cumulants of $X = X_A$ with weights depending on $N, M, \pi, \sigma$, raised to the
$1/n$-th power. While I have not been able to estimate it directly, a trick similar to the
one in Theorem 5.3.4 of \cite{BS98} is applicable. Namely, from Proposition \ref{prop:Poi}
and Corollary \ref{cor:defn} we know that the expressions defining various stochastic
measures do converge to a limit for the free Poisson process. For that process the
expression \eqref{Cauchy} is the sum with all the free cumulants equal to $1$. Since, as
before, in general the cumulants grow no faster than an exponential, we have a uniform
estimate on the expression \eqref{Cauchy}  for $N , M > N_0$ for a fixed $N_0$. Therefore
the sequence of approximations to the stochastic measure is a Cauchy sequence, and so
converges to a limit.

\section{The free Kailath-Segall formula}
\label{sec:KS}

The purpose of this section is to investigate the issues surrounding the free
Kailath-Segall formula, the analog of Theorem 2 of \cite{RW97}.

The free Poisson distribution for the value of the parameter $t := \phi(X_A) = \abs{A} = 1$
is the image of the standard semicircle measure under the squaring map. Therefore if $T_n$
are the orthogonal polynomials with respect to the semicircle measure, namely the Chebyshev
polynomials of the second kind, then the orthogonal polynomials with respect to free
Poisson$(1)$ measure are $P_n (x) = T_{2n}(\sqrt{x})$. In particular from the usual
Chebyshev recursion $x T_n(x) = T_{n+1}(x) + T_{n-1}(x)$ we get the recursion relations
\begin{equation*}
x P_n(x) = P_{n+1}(x) + 2 P_n (x) + P_{n-1}(x).
\end{equation*}
The orthogonal polynomials for the compensated (i.e. centered) Poisson$(1)$ measure satisfy
the same relations. We will see a generalization of these relations in Corollary
\ref{cor:PC}.

Denote $\alpha_{n, m} = \Delta_n \psi_m$ and $\beta_{n, m} = \St_{\hat{1}_n +
\hat{0}_m}$. Let $t := \phi(X)$ be the expectation of the process.
\begin{Lemma}
We have the following recursion relation:
\begin{equation*}
\alpha(n, m) = \beta(n, m) + \sum_{l=0}^{m-1} t^{m-1-l} \beta(n+1, l)
\end{equation*}
for $n, m \geq 1$, with boundary conditions $\alpha(n, 0) = \beta(n, 0)$, $\alpha(0, n) =
\beta(0, n)= \beta(1, n-1)$.
\end{Lemma}

\begin{proof}
The boundary conditions follow directly from the definitions.

\begin{align*}
\Delta_n \psi_m
&=  \lim_{N \rightarrow \infty} \sum_{j=1}^N X_j^n \sum_{\substack{(i_1, i_2, \ldots, i_m) \\
\in \set{1, 2, \ldots, N}^k_{\hat{0}_m}}} X_{i_1} X_{i_2} \cdots X_{i_m} \\
&= \St_{\hat{1}_n + \hat{0}_m} + \lim_{N \rightarrow \infty} \sum_{j=1}^{m}
 \sum_{\substack{(i_1, i_2, \ldots, i_m) \\ \in \set{1, 2, \ldots, N}^m_{\hat{0}_m}}}
 X_{i_j}^n X_{i_1} X_{i_2} \cdots X_{i_m} \\
&= \St_{\hat{1}_n + \hat{0}_m} + \sum_{j=1}^m \St_{\pi_j},
\end{align*}
where $\pi_j = \set{(1, 2, \ldots, n, n + j), (n+1), \ldots, (n+j-1), (n+j+1), \ldots,
(n+m)}$. By an argument similar to the one in the proof of Proposition \ref{prop:Poi},
each inner singleton contributes a factor of $t$, and so $\St_{\pi_j} = t^{j-1}
\St_{\hat{1}_{n+1} + \hat{0}_{m-j}}$.
\end{proof}

\begin{Thm}[free Kailath-Segall formula]
We have the following expression for the $n$-th stochastic measure:
\begin{align*}
\psi_n &= X \psi_{n-1} + \sum_{j=2}^n (-1)^{j-1} \sum_{q=0}^{n-j} \binom{n-q-2}{j-2}
t^{n-j-q} \Delta_j \psi_q \\
 &= X \psi_{n-1} + \sum_{j=2}^n (-1)^{j-1} \sum_{m=0}^{n-j} \binom{m+j-2}{j-2} t^m
\Delta_j \psi_{n-j-m}.
\end{align*}
In particular, $\psi_n$ is a polynomial in the diagonal measures $\Delta_1, \Delta_2,
\ldots, \Delta_n$ and $t$.
\end{Thm}

\begin{proof}
By a repeated use of the Lemma,
\begin{align}
\beta(0, n) &= \beta(1, n-1) = \alpha(1, n-1) - \sum_{l_1 = 0}^{n-2} t^{n-2-l_1}
\beta(2, l_1) \notag \\
 &= \alpha(1, n-1) - \sum_{l_1 = 0}^{n-2} t^{n-2-l_1} \left[ \alpha(2, l_1) -
\sum_{l_2 = 0}^{l_1 - 1} t^{l_1 - 1 - l_2} \beta(3, l_2) \right] \notag \\
 &= \alpha(1, n-1) - \sum_{l_1 = 0}^{n-2} t^{n-2-l_1} \alpha(2, l_1) + \sum_{l_1 = 0}^{n-2}
\sum_{l_2 = 0}^{l_1 - 1} t^{n - 3 - l_2} \beta(3, l_2)  = \ldots \notag \\
 &= \alpha(1, n-1) + \sum_{j=2}^k (-1)^{j-1} \sum_{l_1 = 0}^{n-2} \sum_{l_2 = 0}^{l_1 - 1}
\ldots \sum_{l_{j-1} = 0}^{l_{j-2} - 1} t^{n - j - l_{j-1}} \alpha(j, l_{j-1}) \label{rec} \\
 &\qquad + (-1)^k \sum_{l_1 = 0}^{n-2} \ldots \sum_{l_k = 0}^{l_{k-1} - 1} t^{n - k - 1 - l_k}
\beta(k+1, l_k) \notag.
\end{align}
For $k = n-1$, using the boundary conditions, the last term in \eqref{rec} is
\begin{equation*}
(-1)^k \sum_{l_1 = 0}^{n-2} \ldots \sum_{l_k = 0}^{l_{k-1} - 1} t^{n - (k + 1) - l_k}
\alpha(k+1, l_k).
\end{equation*}
 Thus, continuing the expression \eqref{rec},
\begin{align*}
 &= \alpha(1, n-1) + \sum_{j=2}^n (-1)^{j-1} \sum_{l_1 = 0}^{n-2} \sum_{l_2 = 0}^{l_1 - 1}
\ldots \sum_{l_{j-1} = 0}^{l_{j-2} - 1} t^{n - j - l_{j-1}} \alpha(j, l_{j-1}) \\
 &= \alpha(1, n-1) \\
 &+ \sum_{j=1}^n (-1)^{j-1} \sum_q \abs{\set{l_1, l_2, \ldots, l_{j-1} : q = l_{j-1} < l_{j-2}
< \ldots < l_1 < n-1 }} t^{n-j-q} \alpha(j, p) \\
 &= \alpha(1, n-1) + \sum_{j=1}^n (-1)^{j-1} \sum_{q=0}^{n-j} \binom{n-q-2}{j-2} t^{n-j-q}
\alpha(j, p).
\end{align*}
The result follows from the definitions of $\alpha$ and $\beta$.
\end{proof}

\begin{Cor}
\label{cor:cent}
If the process is centered, i.e. $t=0$, then
\begin{align*}
\psi_n &= \sum_{j=1}^n (-1)^{j-1} \Delta_j \psi_{n-j} \\
 &= \sum_{k=1}^n (-1)^{n-k} \sum_{\substack{j_1, j_2, \ldots, j_k \geq 1 \\ j_1 + j_2 +
\ldots + j_k = n}} \Delta_{j_1} \Delta_{j_2} \cdots \Delta_{j_k}.
\end{align*}
\end{Cor}

\begin{Cor}
For the free Brownian motion X, for $n \geq 2$
\begin{align*}
\psi_n &= X \psi_{n-1} - \abs{A} \psi_{n-2} \\
 &= \sum_{j=0}^{[n/2]} (-1)^j \binom{n-j}{j} \abs{A}^j X^{n-2j}.
\end{align*}
These are the Chebyshev polynomials of the second kind.
\end{Cor}

\begin{Cor}[Recursion relation for the free Poisson process]
For the free Poisson process $X$, for $n \geq 2$,
\begin{equation*}
X \psi_n = \psi_{n+1} + (1-t) \psi_n + t X \psi_{n-1}
\end{equation*}
\end{Cor}

\begin{proof}
\begin{align*}
\psi_n &= X \psi_{n-1} + \sum_{j=2}^n (-1)^{j-1} \sum_{q=0}^{n-j} \binom{n - q - 2}{j - 2}
t^{n-j-q} X \psi_q \\
 &= X \psi_{n-1} + \sum_{q=0}^{n-2} \sum_{j=2}^{n-q} (-1)^{(j-2) + 1} \binom{n-q-2}{j-2}
 t^{(n-q-2)-(j-2)} X \psi_q \\
 &= X \psi_{n-1} - \sum_{q=0}^{n-2} (t-1)^{n-q-2} X \psi_q.
\end{align*}
Therefore
\begin{equation*}
(t-1) \psi_{n-1} = (t-1) X \psi_{n-2} - \sum_{q=0}^{n-3} (t-1)^{n-q-2} X \psi_q
\end{equation*}
and so
\begin{align*}
\psi_n &= (t-1) \psi_{n-1} - X \psi_{n-2} - (t-1) X \psi_{n-2} + X \psi_{n-1} \\
 &= (X + t -1) \psi_{n-1} - t X \psi_{n-2}.
\end{align*}
\end{proof}

Now we consider the compensated Poisson process. From Corollary \ref{cor:poi-d_n} it
follows that the diagonal measures of this process are linear functions in $X$. Therefore
by Corollary \ref{cor:cent}, the stochastic measures $\psi_n$ are $n$-th degree polynomials
in $X$. By Proposition \ref{prop:orth}, they are precisely the polynomials orthogonal with
respect to the compensated free Poisson measure. It is natural to call them free
Poisson-Charlier polynomials. Note that these polynomials have appeared in \cite{HT98}.

\begin{Cor}[Recursion relation for the free Poisson-Charlier polynomials]
\label{cor:PC}
Let $X$ be the free Poisson process, and $\psi_i$ be the stochastic measures for the free
compensated Poisson process $X - t$. Then
\begin{equation*}
X \psi_n = \psi_{n+1} + (1+t) \psi_n + t \psi_{n-1}.
\end{equation*}
\end{Cor}

\begin{proof}
Since the diagonal measures for the free compensated Poisson process are $\Delta_1 = X -
t$, $\Delta_i = X, i \geq 2$, by Corollary \ref{cor:cent}
\begin{equation*}
\psi_n = - t  \psi_{n-1} + \sum_{j=0}^{n-1} (-1)^{n-j-1} X \psi_j.
\end{equation*}
Therefore
\begin{equation*}
\psi_n + \psi_{n-1} = - t \psi_{n-1} - t \psi_{n-2} + X \psi_{n-1}.
\end{equation*}
\end{proof}

\begin{Cor}
The free Poisson-Charlier polynomials are
\begin{equation*}
\psi_n = (X-t)^n + \sum_{i=0}^{n-2} (X-t)^i \sum_{k=1}^{\left[ \frac{n-i}{2}\right]}
\binom{n-i-k-1}{k-1} (-1)^{n-k-i} X^k.
\end{equation*}
\end{Cor}

\begin{proof}
By Corollary \ref{cor:cent},
\begin{align*}
\psi_n &= \sum_{k=1}^n (-1)^{n-k} \sum_{\substack{j_1, j_2, \ldots, j_k \geq 1 \\
j_1 + j_2 + \ldots + j_k = n}} \Delta_{j_1} \Delta_{j_2} \cdots \Delta_{j_k} \\
 &= \sum_{i=0}^n (X-t)^i \sum_k (-1)^{n-k-i} X^k \sum_{\substack{j_1, j_2, \ldots, j_k \geq 2 \\
j_1 + j_2 + \ldots + j_k = n-i}} 1 \\
 &= (X-t)^n + \sum_{i=0}^{n-2} (X-t)^i \sum_{k=1}^{\left[ \frac{n-i}{2}\right]}
\binom{n-i-k-1}{k-1} (-1)^{n-k-i} X^k.
\end{align*}
\end{proof}

\begin{Cor}
For a free compound Poisson process $X$ with generator $e$, $X = s e s$,
\begin{equation*}
\psi_n = X \psi_{n-1} - \sum_{q=0}^{n-2} s (t-e)^{n-q-2} e^2 s \psi_q.
\end{equation*}
\end{Cor}

\begin{proof}
By Corollary \ref{cor:cPoi}, the diagonal measures of $X$ are $\Delta_i = s e^i s$.
Therefore by the Theorem
\begin{align*}
\psi_n &= X \psi_{n-1} + \sum_{j=2}^{n} \sum_{q=0}^{n-j} (-1)^{j-1} \binom{n-q-2}{j-2}
t^{(n-j-q)} s e^j s \psi_q \\
&= X \psi_{n-1} + \sum_{q=0}^{n-2} \sum_{j=2}^{n-q} (-1)^{j-1} \binom{n-q-2}{j-2}
t^{(n-q-2) - (j-2)} s e^{j-2} e^2 s \psi_q \\
 &= X \psi_{n-1} - \sum_{q=0}^{n-2} s (t-e)^{n-q-2} e^2 s \psi_q.
\end{align*}
\end{proof}

\section{Integration of functions}
\label{sec:func}

One of the important points in the proof of the free It\^{o} formula of \cite{BS98} is the
observation that (in the language of this paper) for $X$ the free Brownian motion, for $Z$
centered $\lim_{N \rightarrow \infty} \sum_{i=1}^N X_i Z X_i = 0$. In fact this is true
more generally.

\begin{Thm}
\label{thm:xzx}
Let $Z$ be centered and freely independent from the process $\set{X_I}$. Then
\begin{equation*}
\lim_{N \rightarrow \infty} \sum_{i=1}^N X_i Z X_i = 0.
\end{equation*}
\end{Thm}

\begin{proof}
First, let $X$ be a free Poisson process, $X = s p s$. Then $X_i Z X_i = s (p_i (s Z s)
p_i) s)$. The joint distribution of $Z$ and $sps$ is determined by the condition that $Z$
is freely independent from $sps$. Thus we may assume that in fact $Z$ is freely independent
from $p$. In that case $sZs$ is freely independent from $p_i$, and $\phi(sZs) = \phi(s^2)
\phi(Z) = 0$. Then the result follows from Theorem \ref{thm:eze}.

In general,

\begin{align}
\phi \left( \bigl(\sum_{i=1}^N X_i Z X_i \bigr)^n \right)
 &= \sum_{i_1, i_2, \ldots, i_n = 1}^N \phi(X_{i_1} Z X_{i_1} X_{i_2} Z X_{i_2} \cdots
X_{i_n} Z X_{i_n}) \notag \\
 &= \sum_{\pi \in \Part(n)} \sum_{\substack{(i_1, i_2, \ldots, i_n) \\
\in \set{1, 2, \ldots, N}_\pi^n}} \phi(X_{i_1} Z X_{i_1} X_{i_2} Z X_{i_2} \cdots
X_{i_n} Z X_{i_n}) \notag \\
 &= \sum_{\pi \in \Part(n)} \sum_{\substack{(i_1, i_2, \ldots, i_n) \\
\in \set{1, 2, \ldots, N}_\pi^n}}
\sum_{\sigma \in \NC(n)} \sum_{\substack{\tau \in \NC(2n) \\ \sigma \cup
\tau \in \NC(3n) \\ \tau \leq \pi^2}} R_\sigma(Z, Z, \ldots Z) \notag \\
 & \qquad \qquad \qquad \qquad
\times R_\tau(X_{i_1}, X_{i_1}, X_{i_2}, X_{i_2}, \ldots, X_{i_n}, X_{i_n}) \notag \\
 &= \sum_{\pi \in \Part(n)} \sum_{\substack{(i_1, i_2, \ldots, i_n) \\
\in \set{1, 2, \ldots, N}_\pi^n}}
\sum_{\sigma \in \NC(n)} \sum_{\substack{\tau \in \NC(2n) \\ \sigma \cup
\tau \in \NC(3n) \\ \tau \leq \pi^2}} R_\sigma(Z) N^{-\abs{\tau}}
R_\tau(X) \notag\\
 &= \sum_{\pi \in \Part(n)} (N)_{\abs{\pi}}
\sum_{\sigma \in \NC(n)} \sum_{\substack{\tau \in \NC(2n) \\ \sigma \cup
\tau \in \NC(3n) \\ \tau \leq \pi^2}} R_\sigma(Z, Z, \ldots Z) N^{-\abs{\tau}}
R_\tau(X, X, \ldots, X) \label{xzx}.
\end{align}
Here in the above equation \eqref{xzx}, for partitions $\sigma \in
\NC(n), \tau \in \NC(2n)$, in the partition $\sigma \cup \tau$ we let $\sigma$ act on $\set{2,
5, \ldots, 3n-1}$ while $\tau$ acts on $\set{1,3,4,5, \ldots, 3n-2, 3n}$.

This sum can probably be estimated directly using the properties of Stirling numbers, but
instead we'll use the method of Section \ref{sec:wd}: the result follows from the fact that
the cumulants grow no faster than an exponential, and the fact that the limit exists for
the free Poisson process.

\end{proof}

\begin{Cor}
\label{cor:xzx}
In general, for $Z$ freely independent from the process $\set{X_I}$
\begin{equation*}
\lim_{N \rightarrow \infty} \sum_{i=1}^N X_i Z X_i = \phi(Z) \langle X, X \rangle,
\end{equation*}
where $\langle X, X \rangle = \Delta_2$ is the
quadratic variation of the process.
\end{Cor}

\begin{Prop}
\label{prop:xzxzx}
For $Z_1, Z_2, \ldots, Z_k$ centered and freely independent from the process $\set{X_I}$,
\begin{equation*}
\lim_{N \rightarrow \infty} \sum_{i=1}^N X_i^{m_1} Z_1 X_i^{m_2} Z_2 \cdots X_i^{m_k} Z_k
X_i^{m_{k+1}} = 0.
\end{equation*}
\end{Prop}

\begin{proof}
As in Theorem \ref{thm:xzx},
\begin{align*}
 & \phi \left( \bigl(\sum_{i=1}^N X_i^{m_1} Z_1 X_i^{m_2} \cdots X_i^{m_k} Z_k X_i^{m_{k+1}}
 \bigr)^n \right) \\
 &= \sum_{i_1, i_2, \ldots, i_n = 1}^N \phi(X_{i_i}^{m_1} Z_1 X_{i_1}^{m_2} \cdots
 X_{i_1}^{m_k} Z_k X_{i_1}^{m_{k+1}} \cdots X_{i_n}^{m_1} Z_1 X_{i_n}^{m_2} \cdots
 X_{i_n}^{m_k} Z_k X_{i_n}^{m_{k+1}}) \\
 &= \sum_{\pi \in \Part(n)} (N)_{\abs{\pi}}
\sum_{\sigma \in \NC(nk)} \sum_{\substack{\tau \in
\NC(n \sum_{i=1}^{k+1} m_i) \\ \sigma \cup
\tau \in \NC(n(k + \sum_{i=1}^{k+1} m_i)) \\ \tau \leq \pi^{\sum m_i}}}
R_\sigma(Z_1, \ldots, Z_k, Z_1, \ldots,  Z_k) N^{-\abs{\tau}} R_\tau(X).
\end{align*}
Here in $\sigma \cup \tau$, $\sigma$ acts on the subset $\{m_1 + 1, m_1 + m_2 + 2, \ldots,
\sum_{i=1}^k m_i + k, \sum_{i=1}^{k+1} m_i + k + m_1 + 1, \ldots, (n-1) \sum_{i=1}^{k+1} +
(n-1)k + \sum_{i=1}^k m_i + k\}$ while $\tau$ acts on its complement.

The rest of the proof proceeds as in  Theorem \ref{thm:xzx}.

\end{proof}

\begin{Cor}
\label{cor:xzxzx}
In general, for $Z_1, Z_2, \ldots, Z_k$ freely independent from the process $\set{X_I}$
\begin{align*}
\lim_{N \rightarrow \infty} \sum_{i=1}^N X_i^{m_1} Z_1 X_i^{m_2} Z_2 \cdots X_i^{m_k} Z_k
X_i^{m_{k+1}} &= \phi(Z_1) \phi(Z_2) \cdots \phi(Z_k) \lim_{N \rightarrow \infty}
\sum_{i=1}^N X_i^{m_1 + m_2 + \ldots + m_{k+1}} \\
 &= \phi(Z_1) \phi(Z_2) \cdots \phi(Z_k) \Delta_{\sum_{j=1}^{k+1} m_j}.
\end{align*}
\end{Cor}

These properties should allow us to consider integration with respect to free stochastic
processes. We will return to this subject elsewhere.

\bibliographystyle{amsalpha}

\begin{thebibliography}{VDN92}

\bibitem[Bia97a]{Bia97a}
Philippe Biane, \emph{Free {B}rownian motion, free stochastic calculus and
  random matrices}, Free probability theory (Waterloo, ON, 1995), Fields Inst.
  Commun., vol.~12, Amer. Math. Soc., Providence, RI, 1997, pp.~1--19.

\bibitem[Bia97b]{Bia97c}
Philippe Biane, \emph{Some properties of crossings and partitions}, Discrete
  Math. \textbf{175} (1997), no.~1-3, 41--53.

\bibitem[Bia98]{Bia98}
Philippe Biane, \emph{Processes with free increments}, Math. Z. \textbf{227}
  (1998), no.~1, 143--174.

\bibitem[BS98]{BS98}
Philippe Biane and Roland Speicher, \emph{Stochastic calculus with respect to free
  {B}rownian motion and analysis on {W}igner space},  Probab. Theory Related Fields {\bf 112}
  (1998), no.~3, 373--409.

\bibitem[BLS96]{BLS96}
Marek Bo{\.z}ejko, Michael Leinert, and Roland Speicher, \emph{Convolution and
  limit theorems for conditionally free random variables}, Pacific J. Math.
  \textbf{175} (1996), no.~2, 357--388, \texttt{math.OA/9410054}.

\bibitem[Fag91]{Fag91}
Franco Fagnola, \emph{On quantum stochastic integration with respect to
  ``free'' noises}, Quantum probability \& related topics, QP-PQ, VI, World
  Sci. Publishing, River Edge, NJ, 1991, pp.~285--304.

\bibitem[GSS92]{GSS92}
Peter Glockner, Michael Sch{\"u}rmann, and Roland Speicher, \emph{Realization
  of free white noises}, Arch. Math. (Basel) \textbf{58} (1992), no.~4,
  407--416.

\bibitem[HT98]{HT98}
Uffe Haagerup and Steen Thorbj{\o}rnsen, \emph{Random matrices and {K}-theory
  for exact {$C^\ast$}-algebras}, Odense Preprints No 12, 1998.

\bibitem[Kre72]{Kre72}
G.~Kreweras, \emph{Sur les partitions non crois\'ees d'un cycle}, Discrete
  Math. \textbf{1} (1972), no.~4, 333--350.

\bibitem[KS92]{KS92}
Burkhard K{\"u}mmerer and Roland Speicher, \emph{Stochastic integration on the
  {C}untz algebra ${O}\sb \infty$}, J. Funct. Anal. \textbf{103} (1992), no.~2,
  372--408.

\bibitem[Mar98]{Mar98}
Marcin Marciniak, \emph{On ${Q}$-independence, limit theorems and
  $q$-{G}aussian distribution}, Studia Math. \textbf{129} (1998), no.~2,
  113--135.

\bibitem[Nic95]{Nic95}
Alexandru Nica, \emph{A one-parameter family of transforms, linearizing
  convolution laws for probability distributions}, Comm. Math. Phys.
  \textbf{168} (1995), no.~1, 187--207.

\bibitem[NS96]{NS96a}
Alexandru Nica and Roland Speicher, \emph{On the multiplication of free
  ${N}$-tuples of noncommutative random variables}, Amer. J. Math. \textbf{118}
  (1996), no.~4, 799--837, \texttt{math.OA/9604061}.

\bibitem[NS97]{NS97a}
Alexandru Nica and Roland Speicher, \emph{A ``{F}ourier transform'' for
  multiplicative functions on non-crossing partitions}, J. Algebraic Combin.
  \textbf{6} (1997), no.~2, 141--160.

\bibitem[RW97]{RW97}
Gian-Carlo Rota and Timothy~C. Wallstrom, \emph{Stochastic integrals: a
  combinatorial approach}, Ann. Probab. \textbf{25} (1997), no.~3, 1257--1283.

\bibitem[Spe91]{Spe91}
Roland Speicher, \emph{Stochastic integration on the full {F}ock space with the
  help of a kernel calculus}, Publ. Res. Inst. Math. Sci. \textbf{27} (1991),
  no.~1, 149--184.

\bibitem[Spe98]{Spe98}
Roland Speicher, \emph{Combinatorial theory of the free product with
  amalgamation and operator-valued free probability theory}, Mem. Amer. Math.
  Soc. \textbf{132} (1998), no.~627, x+88.

\bibitem[Voi85]{Voi85}
Dan Voiculescu, \emph{Symmetries of some reduced free product ${C}\sp
  \ast$-algebras}, Operator algebras and their connections with topology and
  ergodic theory (Bu\c steni, 1983), Lecture Notes in Math., vol. 1132,
  Springer, Berlin, 1985, pp.~556--588.

\bibitem[VDN92]{VDN92}
D.~V. Voiculescu, K.~J. Dykema, and A.~Nica, \emph{Free random variables}, CRM
  Monograph Series, vol.~1, American Mathematical Society, Providence, RI,
  1992, A noncommutative probability approach to free products with
  applications to random matrices, operator algebras and harmonic analysis on
  free groups.

\bibitem[vW73]{vWa73}
Wilhelm von Waldenfels, \emph{An approach to the theory of pressure broadening
  of spectral lines}, Probability and information theory, II, Springer, Berlin,
  1973, pp.~19--69. Lecture Notes in Math., Vol. 296.

\end{thebibliography}
\providecommand{\bysame}{\leavevmode\hbox to3em{\hrulefill}\thinspace}

\end{document}